\newif\ifpersonal
\documentclass[12pt, twoside]{article} 
\usepackage{relsize, amsmath,amsthm,mathtools,bm,tensor} 
\usepackage{microtype,lmodern} 
\usepackage[utf8]{inputenc} 
\usepackage[T1]{fontenc} 
\usepackage{enumerate,comment,braket,xspace,tikz-cd,csquotes} 
\usepackage[centering,vscale=0.8,hscale=0.8]{geometry}
\usepackage{xr}
\usepackage[hidelinks,linktoc=all]{hyperref}

\usepackage[garamond,cal=cmcal]{mathdesign}
\usepackage{eucal}

\usepackage[capitalize]{cleveref}
\usepackage{fancyhdr}
\usepackage{lscape}
\usepackage{adjustbox}
\usepackage{faktor}
\numberwithin{equation}{section}
\usepackage{appendix}
\usepackage{stmaryrd}
\usepackage{subfiles}
\theoremstyle{plain}
\newtheorem{thm-intro}{Theorem}[section]
\newtheorem{thm}{Theorem}[section]
\newtheorem{lem}[thm]{Lemma}
\newtheorem{prop}[thm]{Proposition}

\theoremstyle{definition}
\newtheorem{defin}[thm]{Definition}

\newtheorem{eg}[thm]{Example}
\newtheorem{rem}[thm]{Remark}
\newtheorem{construction}[thm]{Construction}

\newtheorem*{note*}{Note}
\makeatletter
\newcommand{\subjclass}[2][2020]{%
	\let\@oldtitle\@title%
	\gdef\@title{\@oldtitle\footnotetext{#1 \emph{Mathematics subject classification.} #2}}%
}
\newcommand{\keywords}[1]{%
	\let\@@oldtitle\@title%
	\gdef\@title{\@@oldtitle\footnotetext{\emph{Key words and phrases.} #1.}}%
}
\makeatother

	\newcommand{\cA}{\mathcal A}

\DeclareMathOperator{\depth}{depth}
\DeclareMathOperator{\linking}{Link}
\DeclareMathOperator{\unzipping}{Unzip}

\newcommand{\rnum}{\mathbb{R}}

\newcommand{\link}[2]{\linking_{{#1}}({#2})}
\newcommand{\unzip}[2]{\unzipping_{{#1}}({#2})}

	\renewcommand{\Im}{\mathrm{Im}}
	\newcommand{\id}{\mathrm{id}}
	\renewcommand{\epsilon}{\varepsilon}
	
\newcommand{\sS}{\mathscr{S}}
\newcommand{\rr}{\mathbb R}
\newcommand{\cI}{{\mathfrak{I}}}
\newcommand{\cT}{{\mathcal{T}}}

\renewcommand{\Im}{\mathrm{Im}}
\renewcommand{\epsilon}{\varepsilon}

\ifpersonal
\newcommand*{\personal}[1]{\textcolor[rgb]{0,0,1}{(Personal: #1)}}
\newcommand*{\todo}[1]{\textcolor{red}{(Todo: #1)}}
\else
\newcommand*{\personal}[1]{\ignorespaces}
\newcommand*{\todo}[1]{\ignorespaces}
\fi
	
	\date{\today}
	
\begin{document}
	\title{Whitney stratifications are conically smooth}
	
	\author{Guglielmo Nocera and Marco Volpe\footnote{The author was supported by the SFB 1085 (Higher Invariants) in Regensburg,
Germany, funded by the DFG (German Science Foundation).
}}
	\maketitle
	\abstract{The notion of conically smooth structure on a stratified space was introduced by Ayala, Francis and Tanaka. This is a very well behaved analogue of a differential structure in the context of stratified topological spaces, satisfying good properties such as the existence of resolutions of singularities and handlebody decompositions. In this paper we prove Ayala, Francis and Tanaka's conjecture that any Whitney stratified space admits a canonical conically smooth structure. We thus establish a connection between the theory of conically smooth spaces and the classical examples of stratified spaces from differential topology.} 

	\tableofcontents

\section{Introduction}

The theory of manifolds with singularities has a long history, that more or less started to consolidate with the definition of Whitney stratified spaces. Some classical references dealing with Whitney stratifications are \cite{thom1969ensembles}, \cite{StratMaps} or \cite{TopStab}, and more recent textbook accounts of the subjects can be found in \cite{Pflaum} and \cite{goresky1988stratified}. Topological spaces admitting a Whitney stratification include algebraic varieties over the real numbers and, more generally, subanalytic subsets of analytic manifolds (see \cite{HiroSub} or \cite{hardt1975topological}). However, Whitney stratified spaces have the disadvantage of not having a good notion of stratified maps between them, which makes it hard to obtain in a clean way desirable results such as some kind of functorial resolution of singularities, or to treat problems of homotopical nature. 
	
To cope with such flaws, in the paper \cite{AFT}, Ayala, Francis and Tanaka have introduced \textit{conically smooth stratified spaces}, a very natural extension of smooth atlases on topological manifolds to the stratified/singular setting. There they also introduce the definition of conically smooth maps, thus defining a category of conically smooth stratified spaces, and prove results such as a functorial resolution of singularities to smooth manifolds with corners. Another important motivation of the authors to develop the theory of conically smooth stratified space was to extend the definition of factorization homology to the setting of manifolds with singularities (see \cite{ayala2017factorization}). Conically smooth stratified spaces can also be used to provide a useful geometric model of the $\infty$-category of $\infty$-categories (see \cite{ayala2019stratified}).

Very roughly, the definition of a conically smooth atlas for a stratified space is built inductively on the \textit{depth} of the stratification (see \cref{defin-depth}), so that the case of depth $0$ corresponds to the usual atlas for smooth manifolds. More generally, it consists of an open covering by stratified spaces of the type $\rnum^n\times C(Z)$, where $Z$ is a compact stratified space with smaller depth equipped with a conically smooth atlas, satisfying some compatibility conditions. This definition is somewhat technical and it needs some time to be developed in full detail, but we will try to review it briefly later in the paper for the reader's convenience.

Although the theory presented in \cite{AFT} is extremely pleasing from a purely formal point of view, it is not really prodigal of examples. For this reason, the authors have conjectured in \cite[Conjecture 1.5.3]{AFT} that any Whitney stratified space should admit a conically smooth atlas. The goal of our paper is to provide a proof of the aforementioned conjecture. It will be convenient to formulate our result in terms of \textit{abstract stratified sets}, introduced by Mather in \cite{TopStab}, as it will make some intermediate steps of our proof much simpler. In short, these are stratified spaces carrying a structure that axiomatizes the construction of a \textit{control datum} for a Whitney stratified space, i.e. a choice of a compatible family of tubular neighbourhoods of strata. The statement of our main result is then as follows.

\begin{thm}\label{thm_in_intro}
    Any $(M,\sS, \mathfrak I)$ abstract stratified set of finite dimension admits a conically smooth atlas.
\end{thm}

Let us give a brief overview of the strategy we adopt to prove \cref{thm_in_intro}. By applying Thom's first isotopy lemma, one may show that $M$ admits an open covering by stratified spaces of the type $\rnum^n\times C(Z)$, where $Z$ is compact and equipped with a structure of an abstract stratified set (see \cref{fiber-Whitney}). The proof of this result goes back to Mather (see \cite[Theorem 8.3]{StratMaps}), but we review it and adapt it to the setting of abstract stratified sets in \cref{Whitney-are-conical}. We refer to the aforementioned covering as the \textit{Thom-Mather charts} of $M$ (see \cref{Thom-Mather-charts}). We then show, by induction on the depth of $(M,\sS, \mathfrak I)$, that the Thom-Mather charts induce a conically smooth atlas. The base of the induction follows by observing that an abstract stratified set of depth $0$ is essentially a smooth manifold. By induction, we may then endow any Thom-Mather chart $\rnum^n\times C(Z)$ with a conically smooth atlas, since the depth of $Z$ is smaller than the one of the ambient space, and therefore one is left to check that the charts are compatible. We then reduce the problem of checking such compatibilities only to two distinct cases. The first one is the case where one has two charts both coming from a tubular neighbourhood of a stratum of maximal depth $X$: this is treated by carefully manipulating the transition maps coming from the smooth structure of $X$, and the distance functions associated to the tubular neighbourhood. The second is the case where the first chart comes as before from $X$, and the second from a stratum $Y$ which contains $X$ in its closure: we treat this situation essentially by observing that, by removing all the strata of maximal depth from $M$, one obtains an abstract stratified space of smaller depth, to which one may apply the inductive hypothesis.

\subsubsection*{Acknowledgments}We wish to thank Pierre Baumann and Peter Haine for their interest and their extremely fruitful comments on the draft of this paper.

\section{Background on stratified spaces}

In the present section, we recall notations and results from the theory of Whitney stratifications and of conically smooth stratifications.

\subsection{Whitney stratifications (Thom, Mather)}

For our review of the theory of Whitney stratifications, we follow \cite{TopStab}, with minimal changes made in order to connect the classical terminology to the one used in \cite{AFT}.

\begin{defin}\label{smooth-stratification}Let $A$ be a smooth manifold, and $M\subset A$ a subset. A \textbf{smooth stratification} $\sS$ of $M$ is a partition of $M$ into subsets $\{M_\alpha\}_{\alpha\in A}=\sS$, such that each $M_\alpha$ is a smooth submanifold of $A$. More generally, if $A$ is a $C^\mu$-manifold, then a $C^\mu$ stratification of a subset $M$ of $A$ is a partition of $M$ into $C^\mu$-submanifolds of $A$.\end{defin}
\begin{rem}In particular, all strata of a smoothly stratified space $M\subset A$ are locally closed subspaces of $M$.\end{rem}

\begin{defin}[Whitney's Condition B in $\rnum^n$]Let $X$ and $Y$ be smooth submanifolds of $\rnum^n$, and let $y\in Y$ be a point. The pair $(X,Y)$ is said to satisfy \textbf{Whitney's Condition B} at $y$ if the following holds. Let $(x_i)\subset X$ be a sequence converging to $y$, and $(y_i)\subset Y$ be another sequence converging to $y$. Suppose that $T_{x_i}X$ converges to some vector space $\tau$ in the $r$-Grassmannian of $\rnum^n$ and that the lines $x_iy_i$ converge to some line $l$ in the $1$-Grassmannian (projective space) of $\rnum^n$. Then $l\subset \tau$.\end{defin}

\begin{defin}[Whitney's condition B]\label{Whitney-B-general} Let $X$ and $Y$ be smooth submanifolds of a smooth $n$-dimensional manifold $A$, and $y\in Y$. The pair $(X,Y)$ is said to satisfy \textbf{Whitney's condition B} at $y$ if there exist a chart $\phi:U\to \rnum^n$ of $A$ around $y$ such that $(\phi(U\cap X), \phi(U\cap Y))$ satisfies Whitney's Condition B at $\phi(y)$.\end{defin}

Recall (e.g. from \cite[Lemma 1.4.4]{Pflaum}) that Whitney's Condition B is invariant under change of charts, and therefore \cref{Whitney-B-general} is well-posed.

\begin{defin}[Whitney stratification]\label{Whitney}Let $A$ be a smooth manifold. A smooth stratification $(M,\sS)$ on a subset $M$ of $A$ is said to be a \textbf{Whitney stratified space} if the following conditions hold \begin{itemize}
		\item (local finiteness) each point has a neighbourhood intersecting only a finite number of strata;
		\item (condition of the frontier) if $Y$ is a stratum of $\sS$, consider its closure $\overline Y$ in $A$. Then we require that $(\overline Y \setminus Y)\cap M$ is a union of strata, or equivalently that if $S\in \sS$ and $S\cap \overline Y\neq \varnothing$, then $S\subset \overline Y$;
		\item (Whitney's condition B) Any pair of strata of $\sS$ satisfies Whitney's condition B when seen as smooth submanifolds of $A$.
\end{itemize}\end{defin}

Given two strata of a Whitney stratification $X$ and $Y$, we say that $X < Y$ if $X\subset \overline Y$. This is a partial order on $\sS$.

A feature of Whitney stratified spaces, perhaps not very evident from the definition, is the existence of compatible tubular neighbourhoods around strata, in a sense that we will now recall. We refer to \cite[Section 6]{TopStab} for more details.

\begin{defin}\label{tubular-nbd}Let $A$ be a manifold and $X\subset A$ be a submanifold. A \textbf{tubular neighborhood} $T$ of $X$ in $A$ is a triple $(E, \varepsilon, \phi)$, where $\pi : E \to X$ is a vector bundle with an inner product $\langle,\rangle$, $\varepsilon: X\rightarrow\rnum_{>0}$ is a smooth function, and $\phi$ is a diffeomorphism of $B_\epsilon=\{e\in E\mid \langle e,e\rangle<\varepsilon(\pi(e))\}$ onto an open subset of $A$, which commutes with the zero section $\zeta$ of $E$:
	$$\begin{tikzcd}
		B_\varepsilon \arrow[rd,"\phi"]&\\
	 X\arrow[u,"\zeta"]\arrow[r,hook]&A.\end{tikzcd}$$\end{defin}
Let $A$ be any smooth manifold, and let $M\subseteq A$ be a subspace equipped with a smooth stratification $\sS$. By \cite[Corollary 6.4]{StratMaps} we see that any stratum of $(M,\sS)$ has a tubular neighbourhood in $A$, but in general the inner products or the projections of tubular neighbourhoods associated to different strata will not interact in any meaningful way. Requiring $(M, \sS)$ to be a Whitney stratification will ensure that certain compatibility conditions for such a tubular neighborhood hold. These conditions are axiomatised in the definition of an \textbf{abstract stratified set} (see \cite[Definition 8.1]{TopStab}). 

\begin{defin}[Abstract stratified set]\label{Abs}
    An \textbf{abstract stratified set} is a triple $(M, \sS, \cI)$ satisfying the following axioms.
    \begin{itemize}
        \item $M$ is a locally compact Hausdorff topological space, having a countable basis for its topology.
        \item $\sS$ is a family of disjoint locally closed subsets of $M$, and $M$ is the union of all its members.
        \item Each element of $\sS$ is a topological manifold equipped with a smooth structure.
        \item $\sS$ is locally finite.
        \item $\sS$ satisfies the condition of the frontier (see the third point in \cref{Whitney}).
        \item $\cI$ is a triple $\{(T_X), (\pi_X), (\rho_X)\}$, where for each $X\in\sS$, $T_X$ is an open neighbourhood of $X$ in $M$, $\pi_X: T_X\rightarrow X$ is a continuous retraction of $T_X$ onto $X$, and $\rho_X: X\rightarrow\rnum_{\geq 0}$ is a continuous function.
        \item $X = \rho_X^{-1}(0)$.
        \item For any stratum $X$ the map
        $$(\pi_{X}, \rho_{X}):T_X\setminus X \rightarrow X\times \rnum_{>0}$$
        is proper.
        \item For any two strata $X$ and $Y$, denote by $T_{X,Y} \coloneqq T_X\cap Y$, $\pi_{X,Y}\coloneqq\pi_{X_{|T_{X,Y}}}$ and $\rho_{X,Y}\coloneqq\rho_{X_{|T_{X,Y}}}$. The map 
        $$(\pi_{X,Y}, \rho_{X,Y}):T_{X,Y}\rightarrow X\times \rnum_{>0}$$
        is a smooth submersion.
        \item For any strata $X$, $Y$ and $Z$, and any $v\in T_{X, Z}\cap T_{Y, Z}\cap \pi_{Y, Z}^{-1}(T_{X, Y})$, we have
        $$\pi_{X, Y}\pi_{Y, Z}(v) = \pi_{X, Z}(v),$$
        $$\rho_{X, Y}\pi_{Y, Z}(v) = \rho_{X, Z}(v).$$
    \end{itemize}
\end{defin}

\begin{defin} Let $(M,\sS,\cI)$ be an abstract stratified set. If $A$ is a smooth manifold and $f : M \to A$ is a continuous mapping, we will say that $f$ is a \textbf{controlled submersion} if the following conditions are satisfied.
\begin{itemize}\item $f|_X : X \to A$ is a smooth submersion for each stratum $X$ of $M$.\item For any stratum $X$, there is a neighborhood $T_X'$ of $X$ in $T_X$ such that $f(v) =
f\pi_X(v)$ for all $v\in T_X'$.
\end{itemize}
\end{defin}



We explain the situation in \cref{Abs} with an example.

\begin{eg}\label{planestrat}
	Let $M$ be the real plane $\rnum^2$ and $\sS$ the stratification given by $$X=\{(0,0)\}$$
	$$Y=\{x=0\}\setminus \{(0,0)\}$$
	$$Z=M\setminus \{x=0\}.$$ We take $\rnum^2$ itself as the ambient manifold.
	Then Mather's construction of the tubular neighbourhoods associated to the strata gives a result like in \cref{plane}.
	\begin{figure}\centering\includegraphics[width=5cm]{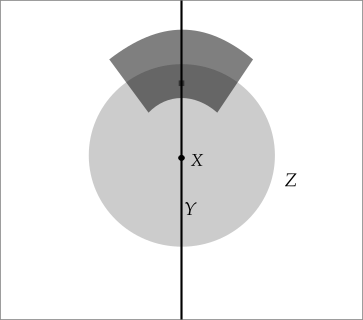}\caption{Tubular neighbourhoods in $(\rnum^2,\{X<Y<Z\})$.}\label{plane}\end{figure}
	Here the circle is $T_X$, and the circular segment is a portion of $T_Y$ around a point of $Y$. We can see here that $T_Y$ is not a ``rectangle'' around the vertical line, as one could imagine at first thought, because the control conditions impose that the distance of a point in $T_Y$ from the origin of the plane is the same as the distance of its ``projection'' to $Y$ from the origin.
	
	Keeping this example in mind (together with its higher-dimensional variants) for the rest of the treatment may be a great help for the visualization of the arguments used in our proofs.
\end{eg}


The following result follows from \cite[Proposition 7.1]{TopStab}.

\begin{prop}\label{Whitney-has-ABS}
Let $(M,\sS)$ be a Whitney stratified space. Then it admits a canonical structure of abstract stratified set in the sense of \cref{Abs}.
\end{prop}

Conversely, we have:

\begin{thm}[{\cite[Theorem on page 3]{Natsume}}]\label{Abs-become-Whitney}
Every paracompact abstract stratified set $(M, \sS, \cI)$ with $\dim V=n$ can be topologically embedded in $\rr^{2n+1}$ such that the image of the embedding is a stratified space satisfying the Whitney condition, and the smooth structures on strata coincide with the ones inherited from $\mathbb R^{2n+1}$. All such ``realizations'' as Whitney stratified spaces in $\mathbb R^N$ are isotopic if $N>2n+2$.
\end{thm}

\begin{rem}
    Notice that, in contrast to \cite[Definition 8.1]{TopStab}, in \cref{Abs} we require the map
    $$(\pi_{X}, \rho_{X}):T_X\setminus X \rightarrow X\times \rnum_{>0}$$
    to be proper. This assumption is verified when $(M, \sS, \cI)$ comes from a Whitney stratified space, and allows one to remove pathological examples such as $C(A)\rightarrow\{0<1\}$ where $A$ is a non-compact manifold. It will also be indispensable for the construction of a conically smooth atlas, as the presence of the latter implies, in particular, that all local links must be compact stratified spaces. 
\end{rem}

\subsection{Conical and conically smooth stratifications (Lurie, Ayala-Francis-Tanaka)}
We now turn to the more recent side of the story, namely the theory of conically smooth stratified spaces. In order to do that, we briefly recall the treatment of stratified sets given by Jacob Lurie in \cite{HA}, which is the base of the formalism used in \cite{AFT}.
\begin{defin}Let $P$ be a partially ordered set. The \textbf{Alexandrov topology} on $P$ is defined as follows. A subset $U \subset P$ is open if it is closed upwards: if $p \leq q$ and $p \in U$ then $q \in U$.\end{defin}
In particular, points are locally closed subsets.
\begin{defin}[{\cite[Definition A.5.1]{HA}}]\label{stratified-HA}A \textbf{stratification} of a topological space $X$ is a continuous map $s:X\to P$ where $P$ is a poset endowed with the Alexandrov topology. The fibers of the points $p\in P$ are subspaces of $X$ and are called the strata. We denote the fiber at $p$ by $X_p$ and by $\sS$ the collection of these strata.\end{defin}

 \begin{rem}Note that, by continuity of $s$, the strata are locally closed subsets of $X$. However, in this definition we do not assume any smooth structure, neither on the ambient space nor on the strata. Furthermore, the condition of the frontier in \cref{Whitney} need not hold for stratifications in the sense of \cref{stratified-HA}: for example, the function $\mathbb R\to \{0<1\}$ given by mapping the interval $(0,1)$ to $1$ and the rest to $0$ is a stratification in the sense of \cref{stratified-HA}, but the two strata do not satisfy the condition of the frontier.

Note however that the condition of the frontier implies that any Whitney stratified space is stratified in the sense of \cref{stratified-HA}: indeed, one obtains a map towards the poset $\sS$ defined by $S<T\iff S\subset \overline T$, which is easily seen to be continuous by the condition of the frontier.\end{rem}
\begin{defin}A \textbf{stratified map} between stratified spaces $(X,P,s)$ and $(Y,Q,t)$ is the datum of a continuous map $f:X\to Y$ and an order-preserving map $\phi:P\to Q$ making the diagram $$\begin{tikzcd}X\arrow[r,"f"]\arrow[d,"s"] & Y\arrow[d,"t"]\\
		P\arrow[r,"\phi"]& Q\end{tikzcd}$$ commute.\end{defin}
\begin{defin}\label{cone-def}Let $(Z,P,s)$ be a stratified topological space. We define \textbf{the cone on $Z$}, denoted by $C(Z)$, as a topological space whose underlying set $$\frac{Z\times[0,1)}{Z\times\{0\}}.$$ Its topology and stratified structure are defined in \cite[Definition A.5.3]{HA}. By \cite[Remark A.5.4]{HA}, when $Z$ is compact and Hausdorff, then the topology is the quotient topology. Note that the stratification of $C(Z)$ is over $P^{\triangleleft}$, the poset obtained by adding a new initial element to $P$: the stratum over this new point is the vertex of the cone, and the other strata are of the form $X\times (0,1)$, where $X$ is a stratum of $Z$.\end{defin}

\begin{defin}[{\cite[Definition A.5.5]{HA}}]\label{Lurie-conical}Let $(X, P, s)$ be a stratified space, $p\in P$, and $x\in X_p$. Let $P_{>p}=\{q\in P\mid q >p\}$. A \textbf{conical chart} at $x$ is the datum of a stratified space $(Z,P_{>p},t)$, an unstratified space $Y$, and a $P$-stratified open embedding $$\begin{tikzcd}Y\times C(Z)\arrow[rr, hook]\arrow[rd]& & X\arrow[ld]\\ &P&\end{tikzcd}$$ whose image hits $x$. Here the stratification of $Y\times C(Z)$ is induced by the stratification of $C(Z)$, namely by the maps $Y\times C(Z)\to C(Z)\to P_{\geq p}\to P$ (see \cref{cone-def}).

A stratified space is \textbf{conically stratified} if it admits a covering by conical charts.\end{defin}

More precisely, the conically stratified spaces we are interested in are the so-called \textbf{$C^0$-stratified spaces} defined in \cite[Definition 2.1.15]{AFT}. Here we recall the two important properties of a $C^0$-stratified space $(X,s:X\to P)$:
\begin{itemize}
	\item every stratum $X_p$ is a \textit{topological manifold};
	\item there is a \textit{basis} of the topology of $X$ formed by conical charts $$\rnum^i\times C(Z)\to X$$ where $Z$ is a \textit{compact} $C^0$-stratified space over the relevant $P_{>p}$. Note that $Z$ will have depth strictly less than $X$; this observation will be useful in order to make many inductive arguments work.
\end{itemize}

Hence the definition of \cite{AFT} may be interpreted as a possible analogue of the notion of topological manifold in the stratified setting: charts are continuous maps which establish a stratified homeomorphism between a small open set of the stratified space and some ``basic'' stratified set.

Following this point of view, one may raise the question of finding an analogue of ``smooth manifold'' (or, more precisely, ``smoothly differentiable structure'') in the stratified setting. We refer to \cite[Definition 3.2.21]{AFT} for the definition of a \textbf{conically smooth structure} (and to the whole Section 3 there for a complete understanding of the notion), which is a very satisfying answer to this question. The definition is rather involved, as it relies essentially on a technical inductive construction based on the \textbf{depth} of a stratification. We sketch the main steps of the definition below.

\begin{defin}[{\cite[Definition 2.4.4]{AFT}}]\label{defin-depth}Let $s : X\rightarrow P$ be a stratified topological space. We define $$\depth(X, P, s) = \sup\limits_{x\in X}(\text{dim}_x(X) - \text{dim}_x(X_{s(x)})),$$ where dim denotes the covering dimension and $X_{s(x)}$ is the stratum of $X$ corresponding to $s(x)\in P$. When the stratifying poset is not specified, and the stratification is given by a family of subsets denoted by $\sS$, we will write $\depth(X, \sS)$. \end{defin}  

\begin{rem}Let $Z$ be an unstratified space of Lebesgue covering dimension $n$. Then the depth of the stratified space $C(Z)\rightarrow\{0<1\}$ at the cone point is $n+1$.\end{rem}

\begin{rem}Note that there is an alternative natural definition of depth at a point $x$: namely, the maximal $k$ such that there exists a chain of strata of the form $X_0<\dots < X_k$ such that $x\in X_0$. The depth of a stratification at a point in the sense of \cref{defin-depth} is always greater or equal than the latter, and they coincide if and only if there is a chain of maximal length $X_0<\dots < X_k$ such that $x \in X_1$ and $\dim X_{i+1}=\dim X_{i}+1$ for $i=0,\dots,k-1$.\end{rem}

The following observation will be useful for the proof of our main result.
\begin{lem}\label{conical-of-depth-zero-is-discrete}
    Let $(X, P,s)$ be a $C^0$ stratified space whose strata are all finite dimensional, and suppose that $\depth(X,P,s)=0$. Then $P$ is a discrete poset, i.e. for any two $p,q\in P$ we have that $p\leq q\iff p=q$.
\end{lem}
\begin{proof}
    Indeed, by assumption, we know that for any $x\in X$ the following formula holds: 
        $$\dim_x X_p =\dim_x X$$
        where $X_p$ is the stratum containing $x$. Up to taking a conical chart centered at $x$, this translates into $$\dim_x(\rr^n)=\dim_{(x,*)}(\rr^n\times C(Z))$$ i.e.
        $$n=n+\dim C(Z)$$
        where $*$ is the cone point of $C(Z)$. This implies that $Z$ is empty, hence the conclusion.
\end{proof} 

   In particular, if $X\to P$ is a $C^0$ stratified space such that $\depth{X} = 0$, then a conically smooth atlas on $X$ is just the usual notion of a smooth atlas which defines $C^{\infty}$-manifolds. If $\depth{X} > 0$, roughly a conically smooth atlas is a collection of stratified open subsets $\{U_i\}_{i\in I}$ of $X$ satisfying the following conditions
   \begin{enumerate}[(i)]
     \item for each $U_i$ there exists a stratified open embedding $\varphi_i : \rnum^{n_i}\times C(Z_i)\hookrightarrow X$ whose image is $U_i$, where $Z_i$ is a compact stratified space equipped with a conically smooth atlas (notice that $\depth{Z_i}<\depth{X}$, so that by induction the notion of an atlas on $Z_i$ is well defined): an object of the type $\rnum^{n_i}\times C(Z_i)$ where $Z_i$ is as above will be called a \textit{basic conically smooth stratified space};
     \item for any $i, j$ such that $U_i\cap U_j\neq\emptyset$, there exists some $k\in I$ and \textit{conically smooth open embeddings} $\rnum^{n_k}\times C(Z_k)\hookrightarrow\rnum^{n_i}\times C(Z_i)$ and $\rnum^{n_k}\times C(Z_k)\hookrightarrow\rnum^{n_j}\times C(Z_j)$ such that the square
     $$
     \begin{tikzcd}
       \rnum^{n_k}\times C(Z_k) \arrow[r, hook] \arrow[d, hook] & \rnum^{n_i}\times C(Z_i) \arrow[d, hook] \\
       \rnum^{n_j}\times C(Z_j) \arrow[r, hook] & X
     \end{tikzcd}
     $$ commutes.
   \end{enumerate}
   For the above definition to make sense, one has to explain what a conically smooth map between basics is. By playing again with inductive arguments on the depth of the target, the main new conceptual input one has to give to formulate precisely this definition is the notion of \textit{differentiability along a cone locus} (see \cite[Definition 3.1.4]{AFT}): for a map $f: \rnum^n\times C(Z)\rightarrow\rnum^m\times C(S)$ between basics, this amounts to requiring the existence of a continuous extension
   $$
   \begin{tikzcd}[column sep  = huge]
     \rnum_{\geq 0}\times \text{T}\rnum^n\times C(Z) \arrow[r, dashed] & \rnum_{\geq 0}\times \text{T}\rnum^m\times C(S) \\
     \rnum_{> 0}\times \text{T}\rnum^n\times C(Z) \arrow[r, "{\gamma^{-1}\circ f_{\Delta}\circ \gamma}"] \arrow[u] & \rnum_{> 0}\times \text{T}\rnum^m\times C(S) \arrow[u]
   \end{tikzcd}
   $$ 
   where the lower horizontal map is built out of $f$ and an appropriate use of the action of scaling and translating on the conical charts such that, in case in which $Z = S = \emptyset$, this recovers the usual notion of differentiability.
   
One main feature of conically smooth structures is the \textit{unzip} construction, that allows one to functorially resolve any conically smooth stratified space into a manifold with corners. For example, if $X_k\hookrightarrow X$ is the inclusion of a stratum of maximal depth, there is a square
   
   \begin{equation}\label{blowupstrata}
     \begin{tikzcd}
       \link{k}{X} \arrow[d, "{\pi_{X}}"] \arrow[r, hook]  & \unzip{k}{X} \arrow[d] \\
       X_k \arrow[r, hook] & X
     \end{tikzcd}
   \end{equation}
 which is both a pushout and a pullback, and $\unzip{k}{X}$ is a conically smooth manifold \textit{boundary} given by $\link{k}{X}$ such that both its interior and $\link{k}{X}$ have depth strictly smaller than the one of $X$. An interesting consequence of the existence of pushout/pullback square \eqref{blowupstrata} is that the notion of conically smooth map is completely determined by the one of smooth maps between manifolds with corners.

\begin{defin} A $C^0$-stratified space together with a conically smooth structure is called a \textbf{conically smooth stratified space}.\end{defin}
 
 \begin{rem}Let us list some useful properties of conically smooth structures:
\begin{itemize}\item any conically smooth stratified space is a $C^0$-stratified space;
	\item all strata have an induced structure of \textit{smooth} manifold, like in the case of Whitney stratifications;

	\item the definition of conically smooth space is intrinsic, in the sense that it does not depend on a given embedding of the topological space into some smooth manifold, in contrast to the case of Whitney stratifications (see \cref{smooth-stratification} and \cref{Whitney});
	\item the notion of conically smooth map (which is a map inducing conically smooth maps between basics in charts) differs substantially from the ``naive'' requirement of being stratified and smooth along each stratum that one has in the case of Whitney stratifications. The introduction of this notion defines a category $\mathcal{S}\text{trat}$ of conically smooth stratified spaces. In this setting, \cite{AFT} are able to build up a very elegant theory and prove many desirable results such as a functorial resolution of singularities to smooth manifolds with corners and the existence of tubular neighbourhoods of conically smooth submanifolds. These results allow to equip $\mathcal{S}\text{trat}$ with a Kan-enrichment (and hence, a structure of $\infty$-category); also, the hom-Kan complex of conically smooth maps between two conically smooth spaces has the ``correct'' homotopy type (we refer to the introduction to \cite{AFT} for a more detailed and precise discussion on this topic), allowing to define a notion of tangential structure naturally extending the one of a smooth manifold and to give a very simple description of the exit-path $\infty$-category of a conically smooth stratified space.
\end{itemize}
\end{rem}
Up to now, the theory of conically smooth spaces has perhaps been in need of a good quantity of explicit examples, especially of topological nature. The following result (conjectured in {\cite[Conjecture 1.5.3]{AFT}}) goes in the direction of providing a very broad class of examples coming from differential geometry and topology.

\begin{thm}\label{main-conjecture} Let $(M,\sS)$ be a Whitney stratified space. Then it admits a conically smooth structure in the sense of \cite{AFT}.\end{thm}

The rest of the paper is devoted to the proof of this theorem (\cref{Whitney-are-conically-smooth}). 
\section{Whitney stratifications are conically smooth}

This section is devoted to proving that any abstract stratified set $(M, \sS, \cI)$ admits a conically smooth atlas (\cref{Whitney-are-conically-smooth}). We will start by constructing conical charts on $M$, essentially adapting the proof of \cite[Theorem 8.3]{StratMaps} to the setting of abstract stratified sets (see \cref{Thom-Mather-charts}). We will then prove that these charts form a conically smooth atlas in the sense of \cite[Definition 3.2.10]{AFT} for $M$, as conjectured in \cite[Conjecture 1.5.3 (3)]{AFT}. We then show that this process associates equivalent conically smooth atlases to equivalent abstract stratified sets.

\subsection{Whitney stratifications are conical}
\begin{lem}\label{fiber-Whitney}Let $(M, \sS,\cI)$ be an abstract stratified set, $T$ a smooth unstratified manifold, and let $f:M\to T$ be a controlled submersion. Then for every $p\in T$ the fiber of $f$ at $p$ has a natural structure of abstract stratified set inherited from $M$.\end{lem}
\begin{proof}
We may assume that for any stratum, up to shrinking tubular neighbourhoods (in such a way to obtain an equivalent abstract stratified sets), for all $v\in T_X$, we have $f(v)=f\pi_X(v)$.
We want to show that, for any $x\in T$, the family $$\{f^{-1}(x)\cap T_X\}_{X\in\cI},$$ together with the obvious restricted $\rho$ and $\pi$, defines an abstract stratified set structure on $f^{-1}(x)$. 

First of all we show that, for any other stratum $Y$, the map $$(T_X\setminus X)\cap Y\cap f^{-1}(x)\xrightarrow[]{(\pi_X, \rho_X)} (X\cap f^{-1}(x))\times \mathbb R_{> 0}$$ is a proper submersion. We know that $$(T_X\setminus X)\cap Y\xrightarrow[]{(\pi_X, \rho_X)} X\times \mathbb R_{> 0}$$ is a proper submersion. Therefore, it will suffice to prove that the square
\[
\begin{tikzcd}
    (T_X\setminus X)\cap Y\cap f^{-1}(x)\arrow[r, "{(\pi_X, \rho_X)}"] \arrow[d, hook] & (X\cap f^{-1}(x))\times \mathbb R_{> 0} \arrow[d, hook] \\
    (T_X\setminus X)\cap Y\arrow[r, "{(\pi_X, \rho_X)}"] &  X\times \mathbb R_{> 0}
\end{tikzcd}
\]
is a pullback. We already know that 
\[
\begin{tikzcd}
(X\cap f^{-1}(x))\times\rnum_{>0} \arrow[d] \arrow[r] & X\cap f^{-1}(x) \arrow[d] \\
X\times\rnum_{>0} \arrow[r]                           & X                        
\end{tikzcd}
\]
is a pullback, and therefore it suffices to show that the square
\[\begin{tikzcd}
	{T_X\cap f^{-1}(x)} & {X\cap f^{-1}(x)} \\
	{T_X} & X
	\arrow[from=1-1, to=2-1]
	\arrow[from=1-1, to=1-2]
	\arrow[from=1-2, to=2-2]
	\arrow[from=2-1, to=2-2]
\end{tikzcd}\]
is a pullback. But the pullback is, by definition, $$\{v\in T_X\mid f\pi(v)=x\}$$ which is equal to $$\{v\in T_X\mid f(v)=x\}=T_X\cap f^{-1}(x).$$
Moreover, since $$(\pi_{X}, \rho_{X}):T_X\setminus X \rightarrow X\times \rnum_{>0}$$ is proper and $(T_X\setminus X)\cap f^{-1}(x)$ is closed in $T_X\setminus X$, we get that the restriction 
$$(T_X\setminus X)\cap f^{-1}(x) \rightarrow (X\cap f^{-1}(x))\times \rnum_{>0}$$
is proper, and thus the proof is concluded.
\end{proof}

\begin{lem}\label{open-are-Whitney}Let $(M,\sS,\cI)$ be an abstract stratified set, and $U\subset M$ an open subset. Then $U$ inherits a natural structure of abstract stratified set obtained by intersecting the elements of $\cI$ with $U$.\end{lem}
\begin{proof}This follows from the immediate observation that an open embedding of smooth manifolds is a submersion, applied to the open embeddings $Y\cap U\subset Y$ for every stratum $Y$.\end{proof}

We will now explain how to provide abstract stratified sets with conical charts. Our construction essentially follows the proof of \cite[Theorem 8.3]{StratMaps}, but adapted to the context of abstract stratified sets (whereas the original proof is given for Whitney stratified spaces). This review is also useful to fix some notations.

\begin{lem}[{Thom's first isotopy lemma, \cite[Corollary 10.2]{TopStab}}]\label{submersions} Let $(M,\sS,\cI)$ be an abstract stratified set, $P$ be a manifold, and $f : M \to P$ be a proper, controlled submersion. Then $f$ is a locally trivial fibration.\end{lem}

\begin{construction}\label{Thom-Mather-charts}
    Let $(M,\sS,\cI)$ be an abstract stratified set. Let $W$ be a stratum of $M$, and let $x$ be a point of $W$. Denote by $\pi_W : T_W\rightarrow W$ the projection, $\rho_W : T_W\rightarrow \mathbb{R}_{\geq 0}$ the ``distance'' function. Let $\varepsilon:W\to \mathbb R_{>0}$ be any smooth function such that, for every $x\in W$, we have $\varepsilon(x)\in \rho_W(\pi_W^{-1}(x))$. Let us choose a ``closed tubular subneighboorhood'' of $T_W$ associated to $\varepsilon$, i.e. the subspace of $T_W$ given by $$N = \{x\in T_W\mid \rho_W(x)\leq \varepsilon(\pi_W(x))\}.$$ Let also $$A=\{x\in T_W\mid \rho_W(x)=\epsilon(\pi_W(x))\}$$ and $f=\pi_W|_A:A\to W$. \personal{Qui abbiamo corretto perch\'e non giustificavamo il fatto che $f$ \`e propria.} Note that $f$ is a proper controlled submersion, since $(\pi_W, \rho_W): T_W\to W\times \mathbb R_{>0}$ is (by \cref{Abs}) and $f$ is the pullback of this map along the graph of $\varepsilon$ $\Gamma_{\varepsilon} : W\to W\times \mathbb R_{>0}$.
Hence by \cref{fiber-Whitney} the restriction of the stratification of $M$ to any fiber of $f$ has a natural structure of abstract stratified set. 

Consider the  mapping 
$$g: N\setminus W\rightarrow W\times (0, 1]$$ defined by 
$$g(x) = \left(\pi_W(x), \frac{\rho_W(x)}{\varepsilon(\pi_W(x))}\right).$$ The space $N\setminus W$ inherits from $M$ a structure of abstract stratified set (see \cref{open-are-Whitney}) and, by \cite[Lemma 7.3 and above]{TopStab}, the map $g$ is a proper controlled submersion. Thus, since $A = g^{-1}(W\times\{1\})$, by \cref{submersions} one gets a stratified homeomorphism $h$ (with respect to the stratification induced on $A$, see \cref{fiber-Whitney}) fitting in the commuting triangle 
$$
\begin{tikzcd}
N\setminus W \arrow[rd, "g"'] \arrow[rr, "h"] &                 & {A\times(0, 1]} \arrow[ld, "f\times \text{id}"] \\
                                              & {W\times(0, 1]} &                                                
\end{tikzcd}.
$$ Furthermore, since $W = \rho^{-1}(0)\subseteq N$, $h$ extends to a homeomorphism of pairs $$(N,W)\xrightarrow{(h,\id)} (M(f),W),$$ where $M(f)$ is the mapping cylinder of $f$ (we recall that $f:A\to W$ is the projection $(\pi_W|_A)$).

 Now for any euclidean chart $j : \rr^i\hookrightarrow W$ around $x$, the pullback of $f$ along $j$ becomes a projection $\rr^i\times Z \rightarrow \rr^i$. Note that $Z$ is compact by properness of $f$, and has an induced structure of abstract stratified set being a fiber of $f$, as we have noticed above. Finally, $$M(f)\simeq M(\rr^i\times Z \xrightarrow{\mathrm{pr_1}}\rr^i)\simeq \overline C(Z)\times \rr^i,$$
 where $\overline C(Z)$ denotes the closed cone on $Z$. The open cone $C(f)$ is thus of the form $C(Z)\times \rr^i$, and this provides the sought conical chart around $x$.

 Summarizing, for an abstract stratified set $(M,\sS,\cI)$, the procedure explained above provides a family $$\cT=\{\phi:\rr^i\times C(Z)\hookrightarrow M\}$$ indexed by the choice of
\begin{itemize}
\item $W\in \sS$
\item $\psi:\rr^i\hookrightarrow W\textup{ smooth chart in }W$
\item $\varepsilon:W\to \rr_{>0}$.\end{itemize}
We will say that a chart is ``centered at the stratum $W$'' and ``induced by the smooth chart $\psi$ and the function $\varepsilon$''. The image of such a chart is an open subset of $T_W$ which we denote by $N(W,\psi,\epsilon)$. The inverse of the chart is (up to composing with the map induced by $\psi^{-1}$) exactly the map $h$ appearing in \cref{Whitney-are-conical}.

We already know that the family $\cT$ is a covering of $M$, and it will be referred to as the family of \textbf{Thom-Mather charts} associated to the abstract stratified set $(M,\sS,\cI)$.
\end{construction}

\begin{rem}\label{same-stratum-same-form}If two charts are centered at the same stratum $W$, up to changing the $\psi$'s and the $\epsilon$'s one can assume that their sources are the same basic $\rr^n\times C(Z)$ (for the same $n$ and $Z$).\end{rem}

This is summarized in the following theorem (it is due to Mather based on ideas by Thom; our previous discussion was meant to elaborate some details necessary for the proof in \cref{Main-section}).

\begin{thm}[{\cite[Theorem 8.3]{StratMaps}}]\label{Whitney-are-conical}Let $(M,\sS,\cI)$ be an abstract stratified set. Then $M$ admits a covering by conical charts of the type $\mathbb R^i\times C(Z)$, where $Z$ is a compact topological space endowed with a natural structure of abstract stratified set.\end{thm}

\begin{eg}
   As a consequence of \cref{Whitney-are-conical}, by \cite[Definition 2.1.15, Axiom (5)]{AFT}, every abstract stratified set is $C^0$-stratified. 
\end{eg}

\subsection{Whitney stratifications are conically smooth}\label{Main-section}

\begin{thm}[Main Theorem]\label{Whitney-are-conically-smooth} Let $(M,\sS, \mathfrak I)$ be an abstract stratified set, and assume that $M$ has finite Lebesgue covering dimension. Then the Thom-Mather charts form a conically smooth atlas on $(M,\sS)$.\end{thm}

\begin{proof} Since $M$ has finite dimension, the depth of $(M, \sS)$ must be finite, so denote $d = \depth(M,\sS)$. The proof will therefore proceed by induction on $d$ (see \cref{defin-depth}). 

In the case $d=0$, by \cref{conical-of-depth-zero-is-discrete} we know that $M$ is just a disconnected union of strata which are smooth manifolds. Therefore, the claim follows easily, since the family of Thom-Mather charts (\cref{Thom-Mather-charts}) reduces to a collection of smooth atlases for each stratum.

Thus, we may assume by induction that for any abstract stratified set $(M',\sS',\cI')$ with \[\depth(M',\sS')< d\] the Thom-Mather charts form a conically smooth atlas on $(M',\sS')$.
	
Fix a Thom-Mather chart $\phi:\rr^i\times C(Z)\hookrightarrow M$, and denote by $\sS_Z$ the stratification on $Z$. Then by the proof of \cref{Whitney-are-conical} we know that $\depth(Z, \sS_Z)< d$ when considering $Z$ with its induced structure of abstract stratified set. Thus, by the inductive hypothesis, the Thom-Mather charts on $Z$ form a conically smooth atlas, and this implies that the $\mathbb R^i\times C(Z)$ is a basic in the sense of \cite[Definition 3.2.4]{AFT}.

 Hence it remains to prove that the ``atlas'' axiom is satisfied: that is, if $m\in M$ is a point, $u:\rnum^i\times C(Z)\to M$ and $v:\rnum^j\times C(W)\to M$ are Thom-Mather charts with images $U$ and $V$, such that $m\in U\cap V$, then there exist a basic $\rr^k\times C(T)$ and a commuting diagram
	\begin{equation}\label{atlas}\begin{tikzcd}\rr^{n_3}\times C(Z_3)\arrow[r,"f_1"]\arrow[d,"f_2"] & \rr^{n_1}\times C(Z_1)\arrow[d,"\phi_1"]\\
\rr^{n_2}\times C(Z_2)\arrow[r,"\phi_2"] & M\end{tikzcd}\end{equation}
such that $x\in \Im (\phi_1f_1)=\Im(\phi_2f_2)$ and that $f_1$ and $f_2$ are maps of basics in the sense of \cite[Definition 3.2.4]{AFT}.

Let $M_d\subseteq M$ be the union of all strata $X$ of depth exactly $d$ (i.e. such that $\sup_{x\in X} \depth_x(M,\sS)=d$). Since $M\setminus M_d$ has depth strictly less than $d$ and it is open in $M$, by the inductive hypothesis we know that the induced structure of abstract stratified set on $M\setminus M_d$ satisfies the claim (i.e. the Thom-Mather charts form a conically smooth atlas). Therefore, it will suffice to examine the following two cases:
\begin{enumerate}\item[(1)] the two charts are both centered at strata of depth less than $d$

\item[(2)]
the two charts are centered at the same stratum $W$ of depth $d$ (and $m$ may or may not be contained in $W$ or belong to $M_d$)

\item[(3)] one chart is centered at a stratum $Y$ of depth $<d$ and the other is centered at a stratum $X$ of depth $d$ and lying in the closure of $Y$. (One may use \cref{planestrat} as a guiding example, with $m$ a point on $\{x=0\}\setminus\{(0,0)\}$.)
\end{enumerate}

\textbf{Case (1).} Indeed, up to choosing smaller $\varepsilon$'s we can assume that both charts lie in $M\setminus M_d$, and hence apply the inductive hypothesis.

\textbf{Case (2).} We define the new chart as follows. Suppose that the chart $\phi_1:\rr^n\times C(Z)\to M$ is centered at $W$ and induced by the smooth chart $\psi_1:\mathbb R^n\to W$ and the function $\epsilon_1:W\to \rr_{>0}$. Analogously, suppose that $\phi_2:\rr^n\times C(Z)\to M$ is centered at $W$ and induced by the smooth chart $\psi_2:\mathbb R^n\to W$ and the function $\epsilon_1:W\to \rr_{>0}$. Note that we can suppose that the basic has the same form in both cases by \cref{same-stratum-same-form}. Choose a smooth chart $\psi_3$ for $X$ and a smooth function $\epsilon_3:W\to \rr_{>0}$ such that \begin{itemize}\item $\Im(\psi_3)\subset \Im(\psi_1)\cap\Im(\psi_2)\subset W$
\item $\epsilon_3(w)\leq \min(\epsilon_1(w),\epsilon_2(w))$ for any $w\in W$
\item the image of the Thom-Mather chart $\phi_3$ associated to $(W, \psi_3, \epsilon_3)$ contains $m$.
\end{itemize}

We call $i_1,i_2$ the transition maps (coming from the smooth structure of $X$) fitting into the diagram

\[\begin{tikzcd}
	{\mathbb R^n} & {\mathbb R^n} \\
	{\mathbb R^n} & W.
	\arrow["{i_2}"', hook, from=1-1, to=2-1]
	\arrow["{i_1}", hook, from=1-1, to=1-2]
	\arrow["{\psi_1}", hook, from=1-2, to=2-2]
	\arrow["{\psi_2}"', hook, from=2-1, to=2-2]
	\arrow["{\psi_3}", hook, from=1-1, to=2-2]
\end{tikzcd}\]
Let us define the following two maps:
$$f_1:\rr^n\times C(Z)\to \rr^n\times C(Z)$$
$$(v,t,z)\mapsto \left(i_1(v), \frac{\epsilon_3(\psi_1i_1(v))}{\epsilon_1(\psi_1i_1(v))}t,z\right)$$

$$f_2:\rr^n\times C(Z)\to \rr^n\times C(Z)$$
$$(v,t,z)\mapsto \left(i_2(v), \frac{\epsilon_3(\psi_2i_2(v))}{\epsilon_2(\psi_2i_2(v))}t,z\right).$$
\textbf{Claim.} These are maps of basics.

\textit{Proof of Claim.} Let us call, for $j=1,2$, $\eta_j=\frac{\epsilon_3\psi_ji_j}{\epsilon_j\psi_ji_j}:\mathbb R^n\to \mathbb R_{>0}.$ This is a smooth function. For $j=1,2$, the verification of the condition that $f_j$ is a map of basics amounts to check that: 

\begin{itemize}\item $f_j$ is conically smooth along $\rr^n$. That is, that the map $$\rr_{>0}\times T\rr^n\times C(Z)\to\rr_{>0}\times T\rr^n\times C(Z) $$
$$(t,v,x,[s,z])\mapsto \left(t,\frac{i_j(tv+x)-i_j(x)}{t},i_j(x),[\eta_j(x)s,z]\right)$$
extends to $t=0$. This follows from smoothness of $i_j$.
\item The differential $$T_x\rr^n\times C(Z)\to T_{i_j(x)}\rr^n\times C(Z)$$
$$(v,[s,z])\mapsto (D_{v}i_j(x),[\eta_j(x)s,z])$$ is injective. This follows from the fact that $i_j$ is a smooth open embedding.
\item the pullback of the atlas $\cA$ that we are considering on $C_1=\rr^n\times \rr_{>0}\times Z$ to $C_2=\rr^n\times \rr_{>0}\times Z$ along $f_j|_{\rr^n\times \rr_{>0}\times Z}=i_j\times \eta\cdot \id\times \id$ coincides with $\cA$. Indeed, for any chart $c_1=j_1\times l_1\times w_1:\rr^n\times \rr_{>0}\times W_1\to C_1$ there exists a chart $c_2=i_j\times \eta\cdot \id\times w:\rr^n\times \rr_{>0}\times W_2\to C_2$ and a commutative diagram 

\[\begin{tikzcd}
	{\rr^n\times \rr_{>0}\times W_1} & &{\rr^n\times \rr_{>0}\times W_2} \\
	{\rr^n\times \rr_{>0}\times Z} & &{\rr^n\times \rr_{>0}\times Z.}
	\arrow["{f|_{\rr^n\times \rr_{>0}\times Z}}"', from=2-1, to=2-3]
	\arrow["{c_1}"', from=1-1, to=2-1]
	\arrow["{j_1\times l_1\times \id}", dotted, from=1-1, to=1-3]
	\arrow["{c_2}", dotted, from=1-3, to=2-3]
\end{tikzcd}\]
\end{itemize}
The last thing to do is to check that, with the notations of the proof of \cref{Whitney-are-conical} and of \cref{Thom-Mather-charts}, the front square in the diagram

\[\begin{tikzcd}
	{N(W,\psi_3,\epsilon_3)} & {N(W,\psi_1,\epsilon_1)} \\
	{N(W,\psi_2,\epsilon_2)} & {\rr^n\times C(Z)} & {\rr^n\times C(Z)} \\
	& {\rr^n\times C(Z)} & M.
	\arrow[hook, from=1-1, to=2-1]
	\arrow[hook, from=1-1, to=1-2]
	\arrow["{\phi_1}", hook, from=2-3, to=3-3]
	\arrow["{f_2}"', hook, from=2-2, to=3-2]
	\arrow["{\phi_2}"', hook, from=3-2, to=3-3]
	\arrow["{f_1}", hook, from=2-2, to=2-3]
	\arrow["\sim"', "h_2", from=2-1, to=3-2, sloped]
	\arrow["\sim"', "h_3", from=1-1, to=2-2, sloped]
	\arrow["\sim"', "h_1", from=1-2, to=2-3, sloped]
\end{tikzcd}\]
commutes. This follows from the fact that, by an easy computation, for $j=1,2$ the diagram

\[\begin{tikzcd}
	{N(W,\psi_3,\epsilon_3)} & {\rr^n\times C(Z)} \\
	{N(W,\psi_j,\epsilon_j)} & {\rr^n\times C(Z)}
	\arrow[hook, from=1-1, to=2-1]
	\arrow["\sim"', "h_3", from=1-1, to=1-2]
	\arrow["\sim"', "h_j", from=2-1, to=2-2]
	\arrow["{f_j}"', hook, from=1-2, to=2-2]
\end{tikzcd}\]
commutes. This completes the proof of the Claim.

\textbf{Case (3).} We can assume that the point $m$ lies outside $M_d$, and that the image of $\phi_2$ is contained in $M\setminus M_d$. Therefore we have the following diagram:

\[\begin{tikzcd}
	& {\rr^n\times \rr_{>0}\times Z_1} & {\rr^n\times C(Z_1)} \\
	{\rr^m\times C(Z_2)} & {M\setminus M_d} & M
	\arrow[hook, from=2-1, to=2-2]
	\arrow[hook, from=1-2, to=2-2]
	\arrow[hook, from=1-2, to=1-3]
	\arrow[hook, from=1-3, to=2-3]
	\arrow[hook, from=2-2, to=2-3]
\end{tikzcd}\]
By \cite[Lemma 3.2.9]{AFT} (``basics form a basis for basics'') there exists a map of basiscs $$i:\rr^{n'}\times C(Z'_1)\hookrightarrow \rr^n\times C(Z_1)$$ whose image is contained in $\rr^n\times\rr_{>0}\times Z_1$ and contains $m$. This induces a commutative diagram

\[\begin{tikzcd}
	& {\rr^{n'}\times C(Z'_1)} \\
	& {\rr^n\times \rr_{>0}\times Z} & {\rr^n\times C(Z)} \\
	{} & {M\setminus M_d} & M.
	\arrow[hook, from=2-2, to=3-2]
	\arrow[hook, from=2-2, to=2-3]
	\arrow[hook, "\phi_1", from=2-3, to=3-3]
	\arrow[hook, from=3-2, to=3-3]
	\arrow[hook, from=1-2, to=2-2]
	\arrow[hook, "i", from=1-2, to=2-3]
\end{tikzcd}\]
and therefore a span
\[\begin{tikzcd}
	& {\rr^{n_1'}\times C(Z'_1)} & {} \\
	{\rr^{n_2}\times C(Z_2)} & {M\setminus M_d}.
	\arrow[from=1-2, to=2-2]
	\arrow[from=2-1, to=2-2]
\end{tikzcd}\]
To this span we can apply the inductive hypothesis, since $\depth M\setminus M_d <d$. This yields maps of basics $f_1':\rr^{n_3}\times C(Z_3)\to \rr^{n_1'}\times C(Z_1')$ and $f_2: \rr^{n_3}\times C(Z_3)\to \rr^{n_2}\times C(Z_2)$ fitting into the diagram
\[\begin{tikzcd}
	{\rr^{n_3}\times C(Z_3)} & {\rr^{n_1'}\times C(Z'_1)} & {\rr^{n_1}\times C(Z_1)} \\
	{\rr^{n_2}\times C(Z_2)} & {M\setminus M_d} & M
	\arrow[from=1-2, to=2-2]
	\arrow[from=2-1, to=2-2]
	\arrow[from=2-2, to=2-3]
	\arrow["i", from=1-2, to=1-3]
	\arrow[from=1-3, to=2-3]
	\arrow["{f_2}", dotted, from=1-1, to=2-1]
	\arrow["{f_1'}", dotted, from=1-1, to=1-2]
\end{tikzcd}\]
Finally, we define $f_1=i\circ f_1'$, which is a map of basics because both $i$ and $f_1'$ are.\end{proof}

We conclude our paper by showing that, given a Whitney stratified space, the conically smooth atlas constructed in \cref{Thom-Mather-charts} is canonical, in the following sense. First of all, recall that two structures of abstract stratified sets on the same Whitney stratified space, induced by different choices of systems of tubular neighbourhoods, are equivalent as abstract stratified sets (\cite[Sec. 8]{TopStab}). Then:

\begin{lem}\label{canonical}
    Let $(M, \sS, \cI)$ and $(M, \sS, \cI')$ be two equivalent abstract stratified sets, and let $\mathcal{A}$, $\mathcal{A}'$ be the conically smooth atlases induced respectively by $\cI$ and $\cI'$ via \cref{Thom-Mather-charts}. Then $\cA$ and $\cA'$ are equivalent.
\end{lem}
\begin{proof}
    Recall that, by definition, two abstract stratified $(M,\sS,\cI), (M,\sS, \cI')$ sets are equivalent if every pair of tubular neighbourhoods $(T_X,\rho_X,\pi_X)\in \cI, (T_X',\rho_X',\pi_X')\in \cI'$ around the same stratum $X$ admits a common tubular subneighbourhood $T''_X$ such that $\rho_X|_{T_X''}=\rho_X'|_{T_X''}, \pi_X|_{T_X''}=\pi_X'|_{T_X''}$. Now, let us consider two $(T_X,\rho_X,\pi_X)\in \cI, (T_Y',\rho_Y',\pi_Y')\in \cI'$ which intersect nontrivially, but are centered at potentially different strata $X,Y$. 
    
    Let $m\in T_X\cap T'_Y$, and let $\phi : \rnum^n\times C(Z)\rightarrow T_X$ and $\psi':\rnum^m\times C(S)\rightarrow T'_Y$ be Thom-Mather charts belonging respectively to $\cA$ and $\cA'$, such that $m\in\text{Im}(\phi)\cap\text{Im}(\psi)$. Arguing as in Case (2) in the proof of \cref{Whitney-are-conically-smooth}, by choosing an appropriately small euclidean neighbourhood of $X$ and rescaling the conical coordinate through a smooth function $\epsilon : X\rightarrow\rnum_{>0}$, one may find a chart $\phi' : \rnum^n\times C(Z)\rightarrow T''_X$ in $\cA'$ and a commutative square 
    $$
    \begin{tikzcd}
\rnum^n\times C(Z) \arrow[d, "\phi'"] \arrow[r] & \rnum^n\times C(Z) \arrow[d, "\phi"] \\
T''_X \arrow[r, hook]                  & T_X                         
\end{tikzcd}
    $$
    where the upper horizontal arrow is a map of basics. But now both $\phi'$ and $\psi'$ belong to $\cA'$, and thus one can find a basic $U$ and a commutative diagram
    $$
    \begin{tikzcd}
U \arrow[d, dotted] \arrow[r, dotted]   & \rnum^n\times C(Z) \arrow[r] \arrow[rd, "\phi'"] & \rnum^n\times C(Z) \arrow[d, "\phi"] \\
\rnum^m\times C(S) \arrow[rr, "\psi'"'] &                                                  & M                                   
\end{tikzcd}
    $$
    where the dotted arrows are maps of basics. Hence we see that $\cA\cup\cA'$ is a conically smooth atlas, and thus the proof is concluded.
\end{proof}

\bibliographystyle{alpha}
\bibliography{Whitney.bib}	
\vfill

	\textit{Guglielmo Nocera}: LAGA, Universit\'e Sorbonne Paris Nord, 99 Av. J.-B. Cl\'ement, 93430 Villetaneuse, France. 
	
	Email: \texttt{guglielmo.nocera@gmail.com}
	
	\
	
	\textit{Marco Volpe}: Max-Planck-Institut f\"ur Mathematik, Vivatsgasse 7, 53111 Bonn, Germany.
	
	Email: \texttt{volpe.marco95@gmail.com}

\end{document}